\documentclass[a4paper]{amsart}

\title{Coherent presentations of structure monoids and the
  Higman-Thompson groups}
\author{Jonathan A. Cohen}
\address{Department of Computing\\
  Macquarie University\\
  Sydney NSW 2109\\
  Australia
}
\email{jonathan.cohen@anu.edu.au}
\date{\today}

\usepackage[matrix,arrow,curve]{xy}
\CompileMatrices

\usepackage{verbatim}
\newtheorem{theorem}{Theorem}[section]
\newtheorem{lemma}[theorem]{Lemma}
\newtheorem{definition}[theorem]{Definition}

\newtheorem{corollary}[theorem]{Corollary}

\newtheorem{example}[theorem]{Example}

\newcommand{\F}{\ensuremath{\mathcal{F}}}
\newcommand{\Fr}{\ensuremath{\mathbb{F}}}
\newcommand{\V}{\ensuremath{\mathcal{V}}}
\newcommand{\E}{\ensuremath{\mathcal{E}}}
\newcommand{\Th}{\ensuremath{(\V,\F,\E)}}
\newcommand{\T}{\ensuremath{\mathcal{T}}}
\newcommand{\sub}{\ensuremath{\mathrm{sub}}}
\newcommand{\dom}{\ensuremath{\mathrm{dom}}}
\newcommand{\Struct}{\ensuremath{\mathrm{Struct}}}
\newcommand{\of}{\ensuremath{\cdot}}
\newcommand{\hT}{\ensuremath{\widehat{\T}}}
\newcommand{\hV}{\ensuremath{\widehat{\V}}}
\newcommand{\hF}{\ensuremath{\widehat{\F}}}
\newcommand{\hE}{\ensuremath{\widehat{\E}}}
\newcommand{\hTh}{\ensuremath{(\hV,\hF,\hE)}}
\newcommand{\supp}{\ensuremath{\mathrm{supp}}}
\newcommand{\D}{\ensuremath{\mathcal{D}}}
\newcommand{\tensor}{\ensuremath{\otimes}}
\newcommand{\Sing}{\ensuremath{\mathrm{Sing}}}
\renewcommand{\P}{\ensuremath{\mathbb{P}}}
\newcommand{\lmb}{\ensuremath{\mathrm{lmb}}}
\newcommand{\Cn}{\ensuremath{{\mathbb{C}_n}}}
\newcommand{\SCn}{\ensuremath{{\mathbb{SC}_n}}}
\newcommand{\Ob}{\ensuremath{\mathrm{Ob}}}
\newcommand{\arrow}[3]{
  \xymatrix@1{
    {#1} \ar[r]^{#2} & #3
  }
}
\begin{document}

\begin{abstract}
  Structure monoids and groups are algebraic invariants of equational
  varieties. We show how to construct presentations of these objects
  from coherent categorifications of equational varieties, generalising
  several results of Dehornoy. We subsequently realise the higher
  Thompson groups $F_{n,1}$ and the Higman-Thompson groups $G_{n,1}$
  as structure groups.  We go on to obtain presentations of these groups via 
  coherent categorifications of the varieties of higher-order
  associativity and of higher-order associativity and commutativity,
  respectively. These categorifications generalise Mac Lane's pentagon
  and hexagon conditions for coherently associative and commutative
  bifunctors. 
\end{abstract}

\maketitle

\section{Introduction}\label{sec:introduction}

Thompson's group $F$ is a finitely presented infinite simple group
that appears in a number of guises. For us, the most useful
description is that of Brown \cite{Brown:finiteness}, which casts the
elements of $F$ as pairs of finite binary trees having the same number
of leaves, subject to a certain equivalence relation on the
pairs. This description suggests that $F$ may in fact have something
to do with associativity, with the elements representing pairs of
equivalent terms in some free semigroup. This observation turns out
to be fruitful and Dehornoy \cite{Dehornoy:thompson} has exploited it
in order to realise $F$ as an algebraic invariant of the variety of
semigroups and subsequently to construct a ``geometric'' presentation
of $F$.  In a similar manner, Dehornoy realises Thompson's group $V$
as an algebraic invariant of the variety of commutative semigroups and
constructs a geometric presentation of $V$.

The relations in Dehornoy's presentations consist of two
parts. First, there are the so-called geometric relations, which arise
purely from the fact that a semigroup is, in the first instance, a
magma. The second class of relations arise from the particular
equational structure of the variety at hand. In the case of
$F$, one additional class of relations are added corresponding
to the Stasheff-Mac Lane pentagon \cite{MacLane_natural} and in the
case of $V$, the presentation further contains a class of relations
corresponding to the Mac Lane hexagon, which encodes the essential
interaction between associativity and commutativity.

The first goal of this paper is to place Dehornoy's constructions in a
more general context. More precisely, instead of a set with operations
and equations, we consider a category with functors and natural
isomorphisms. Within this setting, Dehornoy's geometric relations
correspond to functoriality and naturality of the associated categorical
structure.  The second class of relations correspond to so-called
``coherence axioms'', which are a collection of equations making the
free categorical structure equivalent to a preorder.  

Dehornoy's relation of $F$ and $V$ to particular equational varieties
is a special case of a more general construction
\cite{Dehornoy:varieties} associating an inverse monoid to any
balanced equational variety. This monoid is termed the ``structure
monoid'' of the variety. In Section \ref{sec:free}, we begin by
recalling the construction from \cite{Dehornoy:varieties}. We then go
on to describe ``categorifications'' of equational varieties and show
that a coherent categorification of an equational variety gives rise
to a presentation of the associated structure monoid. In certain
favourable situations, the structure monoid turns out to be a group
and we show that the construction of a presentation from a coherent
categorification carries over to this setting.

Higman \cite{Higman:thompson} has shown that Thompson's group $V$ is
in fact part of an infinite family of finitely presented groups
$G_{n,r}$, which are either simple or have a simple subgroup of index
$2$. Brown \cite{Brown:finiteness} subsequently showed that Thompson's
group $F$ fits into a similar infinite family $F_{n,r}$.  We
recall the definitions of $F_{n,1}$ and $G_{n,1}$ in Section
\ref{sec:groups}. In Section \ref{sec:structure}, we show that the
groups $F_{n,1}$ arise as the structure groups of $n$-catalan
algebras, which encode a notion of associativity for an $n$-ary
function symbol. Similarly, we show that the groups $G_{n,1}$ arise as
the structure groups of symmetric $n$-catalan algebras, which contain
an action of the symmetric group $S_n$ on the variables of an $n$-ary
function symbol. 

In Section \ref{sec:catalancat}, we construct a coherent
categorification of $n$-catalan algebras, which we call $n$-catalan
categories. The coherence axioms for $n$-catalan categories directly
generalise the Stasheff-Mac Lane pentagon axiom, with a new class of
axioms appearing when $n \ge 3$. Following from the results of
sections \ref{sec:free} and \ref{sec:structure}, we obtain new
presentations for $F_{n,1}$. In Section \ref{sec:scatalancat}, we
construct symmetric $n$-catalan categories and show that these form a
coherent categorification of symmetric $n$-catalan algebras, thus
obtaining new presentations for $G_{n,1}$. As in the case of
$n$-catalan categories, additional classes of coherence axioms are
required when $n \ge 3$.

Throughout this paper, we read $f \of g$ as ``$f$ followed by $g$''.

\section{Free categories and structure monoids}\label{sec:free}

We begin this section by recalling Dehornoy's construction of an
inverse monoid associated to a balanced equational theory
\cite{Dehornoy:varieties}. Following this, we describe a process for
obtaining a categorical version of an equational theory and a method
of constructing a monoid from such a categorification. Finally, we
link the two constructions together by showing that coherent
categorifications give rise to presentations of structure monoids. We
base our analysis at the level of theories, rather than of equational
varieties. While this is seemingly at odds with Dehornoy's result
\cite{Dehornoy:varieties} that structure monoids are independant of
the particular equational presentation of a variety, differing
presentations of the same variety lead to distinct categorifactions
and thence to distinct presentations of the structure monoid. 

\subsection{Structure monoids associated to equational theories}

 For a graded set of function symbols $\F$ and a set $X$, we denote by 
$\Fr_\F(X)$ the absolutely free term algebra generated by $\F$ on
$X$. An equational theory is a tuple $\Th$, where $\V$ is a set
of variables, $\F$ is a graded set of function symbols and $\E$ is an
equational theory on $\Fr_\F(\V)$. A map $\varphi:\V \to \Fr_\F(\V)$
is called a \emph{substitution} and it extends inductively to an
endomorphism $\Fr_\F(\V) \to \Fr_\F(\V)$. By abuse of notation, we
label this latter map by $\varphi$ as well. We use $[\V,\Fr_\F(\V)]$
to denote the set of all substitutions. For a term $s \in \Fr_\F(\V)$
and a substitution $\varphi \in [\V,\Fr_\F(\V)]$, we use $s^\varphi$
to denote the image of $s$ under $\varphi$. The \emph{support} of a
term $s$ is the set of variables appearing in it. A pair of terms
$(s,t)$ is \emph{balanced} if they have the same support.

\begin{definition}
  Given a balanced pair of terms $(s,t)$ in $\Fr_\F(\V)$, we use
  $\rho_{s,t}$ to denote the partial function $\Fr_\F(\V) \to
  \Fr_\F(\V)$ with graph
  \[
  \{(s^\varphi,t^\varphi) ~|~ \varphi \in [\V,\Fr_\F(\V)]\}.
  \]
\end{definition}

For a balanced pair of terms $(s,t)$, the partial function
$\rho_{s,t}$ is functional since the support of $t$ is a subset of the
support of $s$. The stronger restriction that the pair is balanced is
required since we wish to utilise the inverse partial function
$\rho_{t,s}$ as well.

Given an equational theory $\T := \Th$, we use $[\E]$ to denote the
congruence generated by $\E$ on $\Fr_\F(\V)$ and we use $\Fr_\T(\V)$
to denote the quotient $\Fr_\F(\V)/[\E]$. Similarly, we use $[s]$ to
denote the congruence class of a term $s$ in $\Fr_\T(\V)$. It is clear
that $[u] = [\rho_{s,t}(u)]$ for any balanced
pair of terms $(s,t)$ and any term $u \in \dom(\rho_{s,t})$. However,
the collection of all partial maps $\rho_{s,t}$ for 
$(s,t) \in \E$ is not sufficient to generate $[\E]$, since equations
apply to subterms as well. To this end, we introduce translated
versions of the maps $\rho_{s,t}$, that apply to arbitrary subterms.

A subterm $s$ of a term $t$ is naturally specified by the node where its
root lies in the term tree of $t$, which in turn is completely
specified by the unique path from the root of $t$ to the root of $s$
in the term tree. A path in a term tree may be specified by an
alternating sequence of function symbols and numbers, where the
numbers indicate an argument of a function symbol. More formally, we
have the following situation.

For a graded set $\F := \coprod_n \F_n$, we set 
\[
A_\F := \bigcup_n\bigcup_{F \in \F_n} \{(F,1),\dots,(F,n)\}.
\]
The set of \emph{addresses associated to $\F$} is denoted by $A_\F^*$ and
is the free monoid generated by $A_\F$ under concatenation, with the
unit being the empty string $\lambda$. For a term $t \in \Fr_\F(\V)$
and an address $\alpha \in A_\F^*$, we use $\sub(t,\alpha)$ to denote
the subterm of $t$ at the address $\alpha$. Note that $\sub(t,\alpha)$
only exists if the term tree of $t$ contains the path $\alpha$ and
that $\sub(t,\lambda) = t$.

\begin{example}
  Suppose that $\F := \{F,G\}$, where $F$ is a binary function symbol
  and $G$ is a ternary function symbol. Suppose that $\V$ is a set of
  variables. Then, the term $t := F(w,G(x,y,z))$ is in
  $\Fr_\F(\V)$. The term tree of $t$ is:
  
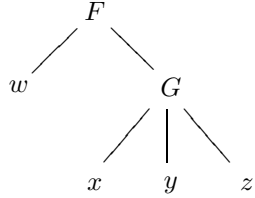
\[
  \begin{xy}
    (0,0)*+{F}="F";
    (-10,-10)*+{w}="w";
    (10,-10)*+{G}="G";
    (0,-23)*+{x}="x";
    (10,-23)*+{y}="y";
    (20,-23)*+{z}="z";
    {\ar@{-} "F"; "w"}
    {\ar@{-} "F"; "G"}
    {\ar@{-} (6,-13); (0,-20)}
    {\ar@{-} (7,-13); (7,-20)}
    {\ar@{-} (8,-13); (14,-20)}
  \end{xy}
  \]
  The term $t$ has the following subterms:
  \begin{center}
    \begin{tabular}{ll}
      \medbreak
      {$\sub(t,(F,1)) = w$} & \quad {$\sub(t,(F,2)) = G(x,y,z)$} \\\medbreak
      {$\sub(t,(F,1)(G,1) = x$} & \quad {$\sub(t,(F,1)(G,2) = y$}\\
      {$\sub(t,(F,1)(G,3)) = z$}
    \end{tabular}
  \end{center}
\end{example} 

\begin{definition}[Orthogonal]
Given a graded set $\F$ and addresses $\alpha,\beta \in A_\F^*$, we
say that $\alpha$ and $\beta$ are orthogonal and write $\alpha \perp
\beta$ if neither $\alpha$ nor $\beta$ is a prefix of the other. Given
a term $t$, and addresses $\alpha$ and $\beta$, the subterms
$\sub(t,\alpha)$ and $\sub(t,\beta)$ are orthogonal if
$\alpha\perp\beta$. 
\end{definition}

Our current addressing system is sufficient to describe translated
copies of the basic operators.

\begin{definition}
  Given a graded set of function symbols $\F$, a variable set $\V$, a
  balanced pair of terms $(s,t) \in \Fr_\F(\V)$ and an address $\alpha
  \in A_\F^*$, the $\alpha$-translated copy of $\rho_{s,t}$ is denoted
  $\rho_{s,t}^\alpha$ and is the partial map $\Fr_F(\V) \to
  \Fr_\F(\V)$ defined as follows:
  \begin{itemize}
  \item A term $u \in \Fr_\F(\V)$ is in the domain of $\rho_{s,t}^\alpha$ if
    $\sub(u,\alpha)$ is defined and is in the domain of $\rho_{s,t}$.
  \item For $u \in \dom(\rho_{s,t}^\alpha)$, the image $\rho_{s,t}^\alpha(u)$
    is defined by \[\sub(\rho_{s,t}^\alpha(u),\alpha) =
    \rho_{s,t}(\sub(u,\alpha))\] and $\sub(\rho_{s,t}^\alpha(u),\beta) =
    \sub(u,\beta)$ for every address $\beta$ orthogonal to $\alpha$.
  \end{itemize}
  Note that $\rho_{s,t}^\lambda = \rho_{s,t}$.
\end{definition}

We are finally in a position to introduce the structure monoid
generated by an equational theory.

\begin{definition}[Structure Monoid]
  Given an equational theory $T := \Th$, the \emph{structure monoid
    of $T$}, denoted $\Struct(T)$, is the monoid of partial
  endomorphisms of $\Fr_\F(\V)$ generated by the following maps under
  composition:
  \[
  \left\{ \rho_{s,t}^\alpha ~|~ (s,t) ~\text{or}~ (t,s) \in \E
    ~\text{and}~ \alpha \in A_\F^* \right\}
  \]
\end{definition}

The structure monoid of an equational theory is readily seen to
completely capture the equational theory.

\begin{lemma}[Dehornoy \cite{Dehornoy:varieties}]\label{lem:capture}
  Let $\T := \Th$ be a balanced equational theory and let $t,t' \in
  \Fr_\F(\V)$. Then $t =_{\T} t'$ if and only if there is some $\rho
  \in \Struct(\T)$ such that $\rho(t) = t'$. \qed
\end{lemma}

Given an equational theory $T = \Th$ and maps $\rho_{s_1,t_1},
\rho_{s_2,t_2} \in \Struct(T)$, the composition
$\rho_{s_1,t_1}\of\rho_{s_2,t_2}$ may be empty. It is nonempty
precisely when there exist substitutions $\varphi,\psi \in
[\V,\Fr_F(\V)]$ such that $t_1^\varphi = s_2^\psi$. In this case, we
say that the pair $(t_1,s_2)$ is \emph{unifiable} and that
$(\varphi,\psi)$ is a \emph{unifier} of the pair. In the case where
$(t_1,s_2)$ is not unifiable, the composition
$\rho_{s_1,t_1}\of\rho{s_2,t_2}$ results in the empty operator, which
we denote by $\varepsilon$. Note that, for any operator $\rho \in
\Struct(T)$, we have $\rho\of\varepsilon = \varepsilon\of\rho =
\varepsilon$. The existence of the empty operator makes freely
computing with inverses in $\Struct(T)$ impossible. 

\begin{definition}[Composable]
  An equational theory $\Th$ is \emph{composable} if any pair of terms in 
  $
  \bigcup_{(s,t)\in\E} \{s,t\}
  $
  are unifiable.
\end{definition}
$\Struct(T)$ always forms an inverse monoid
\cite{Dehornoy:preprint} and contains the empty operator precisely
when $\T$ is not composable. One way in which to transform
$\Struct(G)$ into a group is by passing to the universal group of
$\Struct(\T)$, which we denote by $\Struct_G(\T)$, by collapsing all
idempotents to $1$. In the case where $\T$ is composable, the
idempotent elements of $\Struct(T)$ are precisely those operators that
act as the identity on their domain. A particular class of composable
theories is provided by a certain class of linear theories. Recall
that an equation $s=t$ is \emph{linear} if it is balanced and each
variable appears precisely once in both $s$ and $t$. An equational
theory is linear if each of its defining equations is linear.

\begin{lemma}[Dehornoy \cite{Dehornoy:preprint}]\label{lem:composable}
  A linear equational theory containing precisely one function symbol
  is composable. \qed
\end{lemma}

It follows from the above lemma that each linear equational
theory containing precisely one function symbol gives rise to a
structure \emph{group}.

\begin{example}\label{ex:FV}
 The equational theories for semigroups, $S$, and for commutative
 semigroups, $C$, are both linear. Since these theories involve a single
 binary operator, Lemma \ref{lem:composable} implies that they are
 composable. In this case we have
 that $\Struct_G(S)$ is Thomopson's group $F$ and $\Struct_G(C)$ is
 Thompson's group $V$ \cite{Dehornoy:thompson}. 
\end{example}

\subsection{Categorification of equational theories} An equational
theory $\T = \Th$ defines an algebraic structure on a set. In passing
to the structure monoid $\Struct(T)$ we abstract away from the
underlying set and focus instead on the partial operations generated by
the congruence $[\E]$. This suggests passing to a structure where the
operations generated by $\E$ are given first-class status. In order to
achieve this goal, we pass from a set with algebraic structure to a
category with algebraic structure.

\begin{definition}[Precategorification]
  Given a balanced equational theory $\T = \Th$, the
  \emph{precategorification} of $\T$ is the structure $\hT := \hTh$,
  which consists of:
  \begin{enumerate}
    \item The discrete category $\hV$ generated by $\V$.
    \item For every function symbol $F \in \F_n$, a functor
      $\widehat{F}: \hV^n \to \hV$ in $\hF$. The functor $\widehat{t}$
      for a term $t \in \Fr_\F(\V)$ is defined inductively.
    \item For every equation $(s,t) \in \E$, a natural isomorphism
      $\widehat{\rho}_{s,t}: \widehat{s} \to \widehat{t}$. We use the
      notation $\widehat{\rho}_{t,s} :=
      (\widehat{\rho}_{s,t})^{-1}$.
  \end{enumerate}
\end{definition}

A precategorification of an equational theory should be thought of as
being akin to a graded set of function symbols, rather
than to an equational theory. The reason for this is that, although it
contains all of the information of an equational theory, it does not
contain enough information to ensure that it faithfully represents the
equational structure. A precategorification generates a category whose
objects are the absolutely free term algebra and whose morphisms are
``iterated substitutions'' of the basic maps. Before making this
statement precise, we introduce some notation. Given a term
$\widehat{s} \in \Fr_{\hF}(\hV)$, its support, $\supp(\widehat{s})$, is
the set of objects appearing in it. For a morphism
$\widehat{\rho}:\widehat{s} \to \widehat{t}$, if $\supp(\widehat{s}) =
\supp(\widehat{t}) = \{x_1,\dots,x_n\}$, then we often write
$\widehat{\rho}(x_1,\dots,x_n): \widehat{s}(x_1,\dots,x_n) \to
\widehat{t}(x_1,\dots,x_n)$ to specifically refer to the objects in
the support. Note that a particular $x_i$ may appear more than once in
$\widehat{s}$ or $\widehat{t}$. 

\begin{definition}
  Given a precategorification $\hT := \hTh$ of a balanced equational
  theory $\Th$, we denote by
  $\Fr_{\hT}(\hV)$ the category whose objects are
  $\Fr_{\hF}(\mathrm{Ob}(\hV))$ and whose morphisms are constructed
  inductively as follows:
  \begin{enumerate}
    \item $\Fr_{\hT}(\hV)$ contains the identity $1: \hV \to \hV$.
    \item $\Fr_{\hT}(\hV)$ contains $\hE$.
    \item If $\left\{\widehat{\rho_i} : \widehat{s_i} \to
        \widehat{t_i}\right\}_{i=1}^n \subset \Fr_{\hT}(\hV)$ and $\widehat{F}
      \in \hF_n$, then
      \[
      \widehat{F}(\widehat{\rho_1},\dots,\widehat{\rho_n}) :
      \widehat{F}(\widehat{s_1},\dots,\widehat{s_n}) \to
      \widehat{F}(\widehat{t_1},\dots,\widehat{t_n})
      \]
      is in $\Fr_{\hT}(\hV)$.
    \item If \[\widehat{\rho}(x_1,\dots,x_n): \widehat{s}(x_1,\dots,x_n) \to
      \widehat{t}(x_1,\dots,x_n)\] is in $\Fr_{\hT}(\hV)$ and
      $\widehat{\varphi} \in [\hV,\Fr_{\hF}(\mathrm{Ob}(\hV))]$ is a
      substitution, then
      \[
      \widehat{\rho}(x_1^{\widehat{\varphi}},\dots,x_n^{\widehat{\varphi}}):
      \widehat{s}(x_1^{\widehat{\varphi}},\dots,x_n^{\widehat{\varphi}}) \to 
      \widehat{t}(x_1^{\widehat{\varphi}},\dots,x_n^{\widehat{\varphi}})
      \]
      is in $\Fr_{\hT}(\hV)$
    \item If $\widehat{\rho_1}: \widehat{s} \to \widehat{u}$ and
      $\widehat{\rho_2}:\widehat{u} \to \widehat{t}$ are in
      $\Fr_{\hT}(\hV)$, then
      $\widehat{\rho_1}\of\widehat{\rho_2}:\widehat{s} \to \widehat{t}$
      is in  $\Fr_{\hT}(\hV)$.
  \end{enumerate}
\end{definition}

It is straightforward to show that for $\hT := \hTh$, the category
$\Fr_{\hT}(\hV)$ is the free category on $\hV$ containing all of the
functors in $\hF$ and all of the natural transformations in $\hE$, so
that $\Fr_{\hT}(\hV)$ forms our analogue of the absolutely free term
algebra. 

Categorical structures very rarely arise as precategorifications of
equational theories. Far more common is to require in the definition of
a structure that certain diagrams commute. In particular, given a
collection of diagrams $\D$, each of which consists of a parallel pair
of morphisms $\widehat{\rho_1},\widehat{\rho_2}:\widehat{s} \to
\widehat{t}$  in $\Fr_{\hT}(\hV)$, we may build a congruence $[\D]$ on
the set of morphisms of $\Fr_{\hT}(\hV)$. Factoring out by this
congruence yields the free $\hT$-structure on $\hV$ satisfying the
property that all of the diagrams in $\D$ commute. For particular
categorical structures, it is possible to obtain a set of such
diagrams, whose commutativity implies the commutativity of any diagram
in $\Fr_{(\hT,\D)}(\hV) := \Fr_{\hT}(\hV)/[\D]$. This phenomenon is
termed ``coherence'' and it is equivalent to requiring that
$\Fr_{(\hT,\D)}(\hV)$ is a preorder. Coherence was first
investigated by Mac Lane \cite{MacLane_natural} in relation to
monoidal and symmetric monoidal categories and has subsequently formed
a major part of categorical universal algebra, with an abstract
categorical treatment having been provided by Kelly
\cite{Kelly_coherence}. 

\begin{definition}[Categorification]
  A \emph{categorification} of a balanced equational theory $\T :=
  \Th$ is a pair $(\hT,\D)$, where $\hT$ is a precategorification of
  $\T$ and $\D$ is a collection of parallel pairs of morphisms in
  $\Fr_{\hT}(\hV)$. We say that $(\hT,\D)$ is \emph{coherent} if
  $\Fr_{(\hT,\D)}(\hV)$ is a preorder. 
\end{definition}

\begin{example}\label{ex:monoidal}
  The precategorification of the theory of semigroups consists of a
  binary functor $\tensor$, together with a natural isomorphism:
  \[
  \alpha(x,y,z):x\tensor (y\tensor z) \to (x \tensor y) \tensor z.
  \]
  Mac Lane \cite{MacLane_natural} showed that, in order to obtain a
  coherent categorification of  the theory of semigroups, we need only
  the``pentagon axiom'', which states that the following diagram
  commutes: 
  \[
  \begin{xy}
    (0,0)*+{a\tensor (b\tensor(c\tensor d)))}="1";
    (-25,-18)*+{(a\tensor b)\tensor (c \tensor d)}="2";
    (25,-18)*+{a\tensor((b\tensor c) \tensor d)}="3";
    (-25,-36)*+{((a\tensor b) \tensor c)\tensor d}="4";
    (25,-36)*+{(a\tensor (b\tensor c))\tensor d}="5";
    {\ar@{->}_<<<<<<<{\alpha}@/_/ "1"; "2"}
    {\ar@{->}^<<<<<<<{1\tensor\alpha}@/^/ "1"; "3"}
    {\ar@{->}_{\alpha} "2"; "4"}
    {\ar@{->}^{\alpha} (20,-21); (20,-32)}
    {\ar@{->}^{\alpha\tensor 1} (5,-36); (-12,-36)}
  \end{xy}
  \]
  The precategorification of the theory of commutative semigroups has
  an additional natural isomorphism $\tau$ with the following
  components:
  \[
  \tau(x,y) : x\tensor y \to y\tensor x.
  \]
  Mac Lane \cite{MacLane_natural} went on to show that a coherent
  categorification of the theory of commutative semigroups is obtained
  via the pentagon axiom, together with the axiom that $\tau\of\tau =
  1$ and the ``hexagon axiom'', which states that the following
  diagram commutes:
  \[
  \begin{xy}
    (0,0)*+{a\tensor (b\tensor c)}="1";
    (35,0)*+{(b\tensor c)\tensor a}="2";
    (70,0)*+{b\tensor (c\tensor a)}="3";
    (70,-18)*+{b\tensor (a\tensor c)}="4";
    (0,-18)*+{(a\tensor b)\tensor c}="5";
    (35,-18)*+{(b\tensor a)\tensor c}="6";
    {\ar@{->}^{\tau} "1"; "2"}
    {\ar@{->}^{\alpha^{-1}} "2"; "3"}
    {\ar@{->}^{1\tensor\tau} "3"; "4"}
    {\ar@{->}_{\alpha} "1"; "5"}
    {\ar@{->}_{\tau\tensor 1} (7,-18); (20,-18)}
    {\ar@{->}_{\alpha^{-1}} (40,-18); (53,-18)}
  \end{xy}
  \]
\end{example}

\subsection{Monoid presentations from coherent categorifications}
Dehornoy's utilisation of the pentagon and hexagon coherence axioms in
order to obtain presentations of Thompson's groups
\cite{Dehornoy:thompson} is indicative of a more general relationship
between structure monoids and coherent categorifications of
equational theories. The first step on the road to formalising this
relationship is to construct a monoid presentation out of a
categorification of an equational theory. The initial difficulty is to
construct a set of generators for the monoid corresponding to the
generators of the structure monoid.

\begin{definition}[Singular morphisms]
  Given a balanced equational theory $\T = \Th$ and a
  precategorification $\hT = \hTh$, the set of \emph{singular
    morphisms} of $\Fr_{\hT}(\hV)$ is denoted $\Sing(\hT)$ and is
  generated as follows:  
  \begin{enumerate}
    \item Every natural isomorphism in $\hT$ is singular.
    \item If $\widehat{\rho}$ is singular and $\widehat{\varphi}$ is a
      substitution, then $\hat{\rho}^{\widehat{\varphi}}$ is singular.
    \item If $\widehat{\rho}$ is singular, $\widehat{F} \in
      \widehat{\F_n}$ and $1 \le i \le n$, then
      \[\widehat{F}(\overbrace{1,\dots,1}^{i-1},\widehat{\rho},
      \overbrace{1,\dots,1}^{n-i}),\]
      is singular.
  \end{enumerate}
\end{definition}

In essence, the set of singular morphisms are those that contain
precisely one instance of a generating natural isomorphism. Their
suitability to act as generators is highlighted by the following
lemma.

\begin{lemma}\label{lem:singular}
  Let $\hT := \hTh$ be a precategorification of a balanced equational
  theory $\Th$. Every morphism in $\Fr_{\hT}(\hV)$ is a composite of
  finitely many singular morphisms.
\end{lemma}
\begin{proof}
  The only potential problem is a morphism of the form
  $\widehat{F}(\widehat{\rho_1},\dots,\widehat{\rho_n})$. However, by
  functoriality of $\widehat{F}$, we have 
  \[
   \widehat{F}(\widehat{\rho_1},\dots,\widehat{\rho_n}) = 
   \widehat{F}(\widehat{\rho_1},1,\dots,1)\of\widehat{F}
   (1,\widehat{\rho_2},1,\dots,1)\of\dots\of\widehat{F}
   (1,\dots,1,\widehat{\rho_n}).
  \]
\end{proof}

In order to make the relationship between the monoid we construct from
a categorification and the structure monoid more perspicuous, we
introduce an addressing system for singular morphisms.

\begin{definition}[Type/Address]
  The \emph{type}, $T(\widehat{\rho})$ of a singular morphism
  $\widehat{\rho}$ of a precategorification $\hTh$ is the generating
  natural isomorphism that appears as in instance in it. The address,
  $A(\widehat{\rho})$ is the word of $A^*_{\hF}$ constructed as follows:
  \[
  A(\widehat{\rho}) = \begin{cases}
    (\widehat{F},i)A(\widehat{\sigma}) &\text{if
      $\widehat{\rho}=\widehat{F}(\overbrace{1,\dots,1}^{i-1},\widehat{\sigma},
      1,\dots, 1)$}\\
    \lambda &\text{otherwise}
  \end{cases}
  \]
\end{definition}

Given a categorifacation $(\hT,\D)$, we can now construct a monoid
whose generators are the singular morphisms of $\hT$ and whose
relations are generated by functoriality, naturality and the diagrams
in $\D$. 

\begin{definition}\label{defn:monoid}
  Let $\T := \Th$ be a balanced equational theory and let $(\hT,\D)$
  be a categorification of $\T$. The monoid $\P(\hT,\D)$ is the monoid
  generated by 
  \[
  \{T(\widehat{\rho})^{A(\widehat{\rho})} ~|~ \widehat{\rho} \in
  \Sing(\T)\} \cup
  \{\hat{\rho}_{s,s}^{\widehat{\alpha}}~|~ \text{$(s,t)$ or $(t,s)$ in
    $\hE$ and $\widehat{\alpha} \in A^*_{\hF}$}\} 
  \]
  if $\T$ is composable and by
  \[
  \{T(\widehat{\rho})^{A(\widehat{\rho})} ~|~ \widehat{\rho} \in
  \Sing(\T)\} \cup
  \{\hat{\rho}_{s,s}^{\widehat{\alpha}}~|~ \text{$(s,t)$ or $(t,s)$ in
    $\hE$ and $\widehat{\alpha} \in A^*_{\hF}$}\} \cup
  \{\widehat{\varepsilon}\}
  \]
  otherwise, subject to the following relations.
  \begin{itemize}
    \item \textbf{Inverse:}
      \begin{eqnarray*}
        \hat{\rho}^{\widehat{\alpha}}_{s,t}\of
        \hat{\rho}^{\widehat{\alpha}}_{t,s} &=&
        \hat{\rho}_{s,s}^{\widehat{\alpha}}\\
        \hat{\rho}^{\widehat{\alpha}}_{s,s}\of
        \hat{\rho}^{\widehat{\alpha}}_{s,t}  &=&
        \hat{\rho}^{\widehat{\alpha}}_{s,t}\\
        \hat{\rho}^{\widehat{\alpha}}_{s,t} \of
        \hat{\rho}^{\widehat{\alpha}}_{t,t} &=&
        \hat{\rho}^{\widehat{\alpha}}_{s,t}
      \end{eqnarray*}
      \item \textbf{Composition:} If $t_1$ and $s_2$ are not unifiable
        then 
        \[\hat{\rho}^{\widehat{\alpha}}_{s_1,t_1}\of
        \hat{\rho}^{\widehat{\alpha}}_{s_2,t_2} =
        \widehat{\varepsilon}
        \]
      \item \textbf{Empty operator}:
        \begin{eqnarray*}
          \hat{\rho}^{\widehat{\alpha}}_{s,t} \of
          \widehat{\varepsilon} &=& \widehat{\varepsilon}\\
          \widehat{\varepsilon} \of
          \hat{\rho}^{\widehat{\alpha}}_{s,t}
          &=& \widehat{\varepsilon}
        \end{eqnarray*}
      \item \textbf{Functoriality:} For
        $\widehat{\alpha}\perp\widehat{\beta}$: 
        \[
        \hat{\rho}^{\widehat{\alpha}}_{s,t}\of\hat{\rho}^{\widehat{\beta}}_{u,v}
        = \hat{\rho}^{\widehat{\beta}}_{u,v}\of\hat{\rho}^{\widehat{\alpha}}_{s,t}
        \]
      \item \textbf{Naturality:} Suppose that
        $\hat{\rho}_{s,t}$ is a generator and that
        some variable $x$ appears at addresses
        $\widehat{\beta_1},\dots, \widehat{\beta_p}$ in $\hat{s}$ and
        at addresses $\widehat{\gamma_1},\dots,\widehat{\gamma_q}$ in
        $\hat{t}$. Then, for all addresses
        $\widehat{\alpha},\widehat{\delta}$ and each generator
        $\hat{\rho}_{u,v}$:
        \[
        \hat{\rho}_{s,t}^{\widehat{\alpha}}\of\hat{\rho}_{u,v}^{\widehat{\alpha}
          \widehat{\gamma_1}\widehat{\delta}}\of\dots\of
        \hat{\rho}_{u,v}^{\widehat{\alpha}\widehat{\gamma_q}\widehat{\delta}}
        =
        \hat{\rho}_{u,v}^{\widehat{\alpha}\widehat{\beta_1}\widehat{\delta}}\of\dots 
        \of
        \hat{\rho}_{u,v}^{\widehat{\alpha}\widehat{\beta_p}\widehat{\delta}}
        \of\hat{\rho}_{s,t}^{\widehat{\alpha}}        
        \]
      \item \textbf{Coherence:} For $(\sigma_1\of\dots\of\sigma_p,
        \tau_1,\dots,\tau_q) \in \D$, where each $\sigma_i$ and
        $\tau_j$ is singular, set:
        \[
        T(\sigma_1)^{A(\sigma_1)}\of\dots\of T(\sigma_p)^{A(\sigma_p)}
        =
        T(\tau_1)^{A(\tau_1)}\of\dots\of T(\tau_q)^{A(\tau_q)}
        \]
  \end{itemize}
\end{definition}

The relations for functoriality and naturality in $\P(\hT,\D)$ are
adapted from \cite{Dehornoy:preprint}. The functoriality relation is
precisely the requirement that each operator $\widehat{F} \in \hF$ is
a functor. The naturality condition is, in turn, precisely the
requirement that each $\widehat{\rho} \in \hE$ is a natural
transformation. The rather involved
addressing system in the naturality condition is due to the fact that
the same variable may appear multiple times in different positions on
either side of an equation. For naturality, one needs to apply a map
to each of these instances of the variable simultaneously. We now set
about relating $\P(\hT,\D)$ to $\Struct(\T)$. 

\begin{lemma}
  Let $\T$ be a balanced equational theory and let $(\hT,\D)$ be a
  categorification of $\T$. Then $\P(\hT,\D)$ is an inverse monoid.
\end{lemma}
\begin{proof}
  For nonempty $\hat{\rho} := \hat{\rho}_{s_1,t_1}^{\widehat{\alpha_1}}\of
  \dots \hat{\rho}_{s_k,t_k}^{\widehat{\alpha_k}}$, set
  $\hat{\rho}^{-1} := \hat{\rho}_{t_k,s_k}^{\widehat{\alpha_k}}\of
  \dots \hat{\rho}_{t_1,s_1}^{\widehat{\alpha_k}}$. Then it follows
  from the Inverse relations that
  \begin{eqnarray*}
    \hat{\rho}\of\hat{\rho}^{-1}\of\hat{\rho} &=& \hat{\rho}\\
    \hat{\rho}^{-1}\of\hat{\rho}\of\hat{\rho}^{-1} &=& \hat{\rho}^{-1}
  \end{eqnarray*}
  Since we also have that
  $\widehat{\varepsilon}\of\widehat{\varepsilon}\of\widehat{\varepsilon}
  = \widehat{\varepsilon}$, it follows that $\P(\hT,\D)$ forms an
  inverse monoid.
\end{proof}

\begin{theorem}\label{thm:monoid}
  Let $\T$ be a balanced equational theory and let $(\hT,\D)$ be a
  categorification of $\T$. The following map is an epimorphism of
  inverse monoids and it is an isomorphism if and only if $(\hT,\D)$
  is coherent:
  \begin{eqnarray*}
    \P(\hT,\D) &\stackrel{\Theta}{\longrightarrow}& \Struct(\T)\\
    \hat{\rho}_{s_1,t_1}^{\widehat{\alpha_1}}\of\dots\of
    \hat{\rho}_{s_k,t_k}^{\widehat{\alpha_k}} &\longmapsto&
    \rho^{\alpha_1}_{s_1,t_1}\of\dots\of\rho^{\alpha_k}_{s_k,t_k}
  \end{eqnarray*}
\end{theorem}
\begin{proof}
  By construction, $\Theta$ is a homomorphism of inverse monoids. For
  surjectivity, we need only show that every generator
  $\rho_{s,t}^\alpha \in \Struct(\T)$ corresponds to some singular
  morphism $S(\rho^\alpha_{s,t}) \in \hT$. This singular morphism can be
  constructed recursively as follows:
  \[
  S(\rho^\alpha_{s,t}) = \begin{cases}
    \widehat{F}(\overbrace{1,\dots,1}^{i-1},S(\rho^\beta_{s,t}),1,\dots,1)
    & \text{if $\alpha = (\widehat{F},i)\beta$}\\
    \rho_{s,t} &\text{if $\alpha = \lambda$}
  \end{cases}
  \]
  It remains to show that $\Theta$ is faithful if and only if
  $(\hT,\D)$ is coherent. 

  Suppose that $\Theta$ is faithful and let
  $\widehat{\rho_1},\widehat{\rho_2}$ be a parallel pair of morphisms
  in $\hT$. Then $\Theta(\widehat{\rho_1}) =
  \Theta(\widehat{\rho_2})$, since $\widehat{\rho_1}$ and
  $\widehat{\rho_2}$ have the same source and target. Since $\Theta$ is
  faithful, it follows that $\widehat{\rho_1} = \widehat{\rho_2}$.

  Conversely, suppose that $(\hT,\D)$ is coherent and that
  $\Theta(\widehat{\rho_1}) = \Theta(\widehat{\rho_2})$. Then,
  $\widehat{\rho_1}$ and $\widehat{\rho_2}$ have the same
  source and target. Since $(\hT,\D)$ is coherent, it follows that
  $\widehat{\rho_1} = \widehat{\rho_2}$.
\end{proof}

As in the case of structure monoids, when the theory $\T$ is balanced
and composable, we may construct a group $\P_G(\hT,\D)$ from a
categorification $(\hT,\D)$.

\begin{definition}
  Let $\T$ be a balanced composable equational theory and let
  $(\hT,\D)$ be a categorification of $\T$. The group $\P_G(\hT,\D)$
  is generated by
  \[
  \{T(\widehat{\rho})^{A(\widehat{\rho})} ~|~ \widehat{\rho} \in
  \Sing(\T)\},
  \]
  subject to the functoriality, naturality and coherence relations
  from Definition \ref{defn:monoid}, together with the following
  relation:
  \[
  (\hat{\rho}^{\widehat{\alpha}}_{s,t})^{-1} =   \hat{\rho}^{\widehat{\alpha}}_{t,s}.
  \]
\end{definition}

Following the same line of reasoning as in the proof of Theorem
\ref{thm:monoid}, we obtain the following relationship between
$\P_G(\hT,\D)$ and $\Struct_G(\T)$.

\begin{theorem}\label{thm:group}
  Let $\T$ be a balanced, composable equational theory and let
  $(\hT,\D)$ be a categorification of $\T$. The following map is an
  epimorphism of groups and it is an isomorphism if and only if $(\hT,\D)$
  is coherent:
  \begin{eqnarray*}
    \P_G(\hT,\D) &\stackrel{\Theta}{\longrightarrow}& \Struct_G(\T)\\
    \hat{\rho}_{s_1,t_1}^{\widehat{\alpha_1}}\of\dots\of
    \hat{\rho}_{s_k,t_k}^{\widehat{\alpha_k}} &\longmapsto&
    \rho^{\alpha_1}_{s_1,t_1}\of\dots\of\rho^{\alpha_k}_{s_k,t_k}
  \end{eqnarray*} \qed
\end{theorem}

\begin{example}
  In Example \ref{ex:monoidal}, we obtained a coherent categorification
  of the theory of semigroups, S, consisting of Mac Lane's pentagon
  axiom. It follows from Theorem \ref{thm:group} that we can construct
  a presentation for $\Struct_G(S)$. We saw in Example \ref{ex:FV}
  that $\Struct_G(S)$ is isomorphic to Thompson's group $F$ and we
  thereby obtain a presentation for $F$. Similarly, we obtain a
  presentation of Thompson's group $V$ using the pentagon and hexagon
  coherence axioms from Example \ref{ex:monoidal}. The resulting
  presentations are the same as those constructed by Dehornoy
  \cite{Dehornoy:thompson}.
\end{example}

In the following section, we describe generalisations of Thompson's
groups $F$ and $V$ due to Higman \cite{Higman:thompson} and Brown
\cite{Brown:finiteness}. 

\section{The groups $F_{n,1}$ and $G_{n,1}$}\label{sec:groups} 

In the previous section, we have seen that Thompson's groups $F$ and
$V$ arise as structure groups of certain balanced equational theories
and we have subsequently obtained presentations for these groups via
coherent presentations of their associated categorical theories. In
this section, we introduce generalisations of these groups due to
Brown \cite{Brown:finiteness} and Higman \cite{Higman:thompson}, which
we call $F_{n,1}$ and $G_{n,1}$, respectively. In
the following sections, we shall see how the aforementioned process of
constructing presentations for $F$ and $V$ generalises to this broader
class of groups.

There are several paths to defining the groups $F_{n,1}$ and
$G_{n,1}$, all of which relate to the fact that each of these groups
arises as a subgroup of the automorphism group of a Cantor set. Of
the myriad of definitions available, we choose to follow the
description of Brown \cite{Brown:finiteness}, which utilises certain
equivalence classes of pairs of finite rooted trees.

\begin{definition}[Tree]
  The set of $n$-ary trees is defined inductively as follows:
  \begin{itemize}
    \item The graph consisting solely of a single vertex is an $n$-ary
      tree.
    \item If $T_1,\dots,T_n$ are $n$-ary trees then the
      following is also an $n$-ary tree:
      
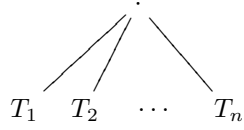
\[
      \begin{xy}
        (0,0)*+{\cdot}="a";
        (-15,-14)*+{T_1}="t1";
        (-7,-14)*+{T_2}="t2";
        (2,-14)*+{\dots}="dots";
        (12,-14)*+{T_n}="tn";
        {\ar@{-} "a"; "t1"};
        {\ar@{-} "a"; "t2"};
        {\ar@{-} "a"; "tn"};
      \end{xy}
      \]
  \end{itemize}
\end{definition}
The \emph{root} of an $n$-ary tree is the unique vertex of valence $0$
or $n-1$. The \emph{leaves} of a rooted tree $T$ are the vertices of
valence $0$ or $1$ and we denote this set by $\ell(T)$. 

\begin{definition}[Expansion]
  A \emph{simple expansion} of an $n$-ary tree $T$ is the tree obtained by
  replacing a leaf $v$ of $T$ with the following:
  
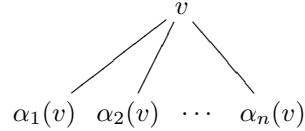
\[
  \begin{xy}
    (0,0)*+{v}="a";
    (-18,-14)*+{\alpha_1(v)}="t1";
    (-7,-14)*+{\alpha_2(v)}="t2";
    (2,-14)*+{\dots}="dots";
    (12,-14)*+{\alpha_n(v)}="tn";
    {\ar@{-} "a"; "t1"};
    {\ar@{-} "a"; "t2"};
    {\ar@{-} "a"; "tn"};
  \end{xy}
  \]
  An \emph{expansion} of an $n$-ary tree is a tree obtained by making
  finitely many succesive simple expansions.    
\end{definition}

Given two trees $T_1$ and $T_2$ having a common expansion $S$, we say
that $S$ is a \emph{minimal} common expansion if any other expansion $S'$ of
$T_1$ and $T_2$ is an expansion of $S$.

\begin{lemma}[Higman \cite{Higman:thompson}]\label{lem:expansion}
  Any two finite $n$-ary trees have a minimal common expansion. \qed
\end{lemma}

The underlying sets of the groups $F_{n,1}$ and $G_{n,1}$ consist of
certain formal expressions called \emph{tree diagrams}.

\begin{definition}[Tree diagram]
  An \emph{$n$-ary tree diagram} is a triple $(T_1,T_2,\sigma)$, where
  $T_1$ and $T_2$ are $n$-ary trees having the same number of leaves
  and $\sigma$ is a bijection $\ell(T_1) \to \ell(T_2)$. 
\end{definition}

As in the case of trees, we may talk about expansions of tree
diagrams.

\begin{definition}
  A \emph{simple expansion} of an $n$-ary tree diagram
  $(T_1,T_2,\sigma)$ is an $n$-ary tree diagram $(T_1',T_2',\sigma')$
  obtained by the following procedure:
  \begin{itemize}
    \item $T_1'$ is a simple expansion of $T_1$ along the leaf $l$. 
    \item $T_2'$ is the simple expansion of $T_2$ along the leaf
      $\sigma(l)$. 
    \item $\sigma'$ is the bijection $\ell(T_1') \to \ell(T_2')$ defined
      by setting $\sigma'(k) = \sigma(k)$ for $k \in
      \ell(T_1)\setminus\{l\}$ and $\sigma'(\alpha_i(l)) =
      \alpha_i(\sigma(l))$.
  \end{itemize}
  An \emph{expansion} of an $n$-ary tree diagram $(T_1,T_2,\sigma)$ is
  any $n$-ary tree diagram obtained by making finitely many succesive
  simple expansions of $(T_1,T_2,\sigma)$.
\end{definition}

Let $\sim$ be the equivalence relation on the set of $n$-ary tree
diagrams obtained by setting $(T_1,T_2,\sigma) \sim (T_1',T_2',\sigma')$
whenever  $(T_1,T_2,\sigma)$ and $(T_1',T_2',\sigma')$ possess a
common expansion. Let $[(T_1,T_2,\sigma)]$ denote the equivalence
class of $(T_1,T_2,\sigma)$ modulo $\sim$. We call
$[(T_1,T_2,\sigma)]$ an \emph{$n$-ary tree symbol}. 

\begin{definition}
For $n\ge 2$, we set $G_{n,1}$ to be the group whose underlying set is
the collection of 
$n$-ary tree symbols, together with the following group structure:
\begin{itemize}
    \item Given two $n$-ary tree symbols $[(T_1,T,\sigma)]$ and
      $[(T',T_2,\sigma']$, it follows from Lemma
      \ref{lem:expansion} that we may assume that $T = T'$. We
      define their product to be
      \[
      [(T_1,T,\sigma)][(T,T_2,\sigma')] = [(T_1,T_2,\sigma\of\sigma')].
      \]
    \item The inverse of $[(T_1,T_2,\sigma)]$ is
      $[(T_2,T_1,\sigma^{-1})]$.
    \item The unit element is $[(T,T,id)]$.
  \end{itemize}
\end{definition}

It follows from the definitions that any $n$-ary tree is an expansion
of the tree consisting solely of a single vertex. Thus, the leaves of
an $n$-ary tree may be seen as a subset of the free monoid on
$\{1,\dots,n\}$. Therefore, we may order the leaves of the tree
lexicographically, which is equivalent to ordering the leaves 
left-to-right when drawn on a page. We say that an $n$-ary tree
symbol $[(T_1,T_2,\sigma)]$ is \emph{order-preserving} if $\sigma$ is
an isomorphism of ordered sets; that is, if $\sigma$ preserves this 
ordering. 

\begin{definition}
  For $n \ge 2$, we set $F_{n,1}$ to be the subgroup of $G_{n,1}$
  consisting of the order-preserving $n$-ary tree symbols.
\end{definition}

The groups $F_{n,1}$ and $G_{n,1}$ generalise Thompson's original
groups $F$ and $V$, since we have $F_{2,1} \cong F$ and $G_{2,1} \cong
V$. They also share several of the interesting properties of $F$ and
$V$ as surveyed in \cite{Scott:tour}. In the following section, we
shall realise $F_{n,1}$ as the structure group of higher-order
associativity and $G_{n,1}$ as the structure group of higher order
associativity and commutativity.  

\section{$F_{n,1}$ and $G_{n,1}$ as structure groups}\label{sec:structure}

Our goal in this section is to realise $F_{n,1}$ and $G_{n,1}$ as
structure groups. Since both of these groups are built using maps
between $n$-ary trees, we take our set of function symbols to be $\F
:= \{\tensor\}$, where $\tensor$ is an $n$-ary function symbol. For a
set of variables $\V$, there is an obvious bijection between
$\Fr_\F(\V)$ and the set of $n$-ary trees whose leaves are labelled by
members of $\V$. We denote the absolutely free term algebra generated
by $\{\tensor\}$ on the set $\V$ by $\Fr_\tensor(\V)$ and we denote
the free monoid generated by $\V$ by $\V^*$.   

Our basic strategy is to first realise $F_{n,1}$ as a
structure group by constructing an equational theory $\E$ such that
$[\E]$ equates any two terms $t_1,t_2 \in \Fr_\F(\V)$ that contain
precisely the same variables in the same order and such that no
variable appears more  
than once in either $t_1$ or $t_2$. In the binary case, there is an
obvious candidate for $\E$: associativity. So, $\E$ needs to be
an analogue of associativity for $n > 2$. Once we have this
realisation of $F_{n,1}$ we need only add the ability to arbitrarily
permute variables in order to obtain a realisation of $G_{n,1}$ as a
structure group.

\subsection{Catalan Algebras and $F_{n,1}$}

Associativity of a binary function symbol is sufficient to establish
that any two bracketings of the same string are equal. The way
in which one establishes this fact is to show that any bracketing of
a string is equal to the left most bracketing. So, for an $n$-ary
function symbol to be associative, we need equations which imply that
any bracketing of a term is equivalent to the left most one. In order
to simplify notation, for integers $i \le j$, we use the symbol $x_i^j$
to denote the list $x_i,x_{i+1},\dots,x_j$. If $i > j$, then $x_i^j$
is the empty list.

\begin{definition}[$n$-Catalan algebras]
  For $n \ge 2$, the theory of $n$-Catalan algebras consists of an
  $n$-ary function symbol $\tensor$ together with the following
  equations, where $0 < i < n$:
  \[
  \tensor(x_1^i,\tensor(x_{i+1}^{i+n}),x_{i+n+1}^{2n-1}) =
  \tensor(x_1^{i-1},\tensor(x_{i}^{i+n-1}),x_{i+n}^{2n-1})
  \]
  We denote the theory of $n$-Catalan algebras by $C_n$. 
\end{definition}
The reason for the name of $n$-catalan algebras is that the set of all
terms having $k$ occurrences of the symbol $\tensor$ and containing
precisely one variable is in bijective 
correspondence with the set of $n$-ary trees having $k$ internal
nodes, which has cardinality equal to the generalised Catalan number 
$\frac{1}{(n-1)k+1}{nk \choose k}$, \cite{Stanley:enumerative}.
The rather opaque equational theory of $n$-Catalan algebras is
rendered somewhat more understandable by viewing the induced equations
on the term trees which, for $n=3$, yields the following:

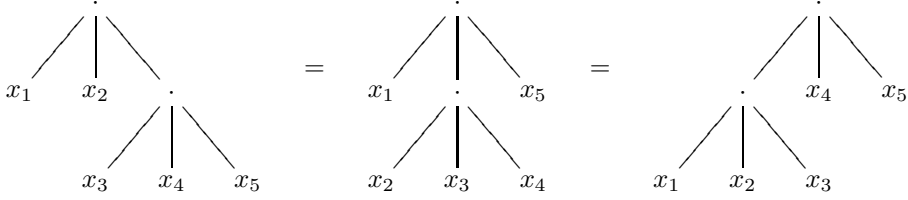
\begin{center}
  \begin{tabular}{clcrc}
    $
    \begin{xy}
      (0,0)*+{\cdot}="a";
      (-10,-12)*+{x_1}="x1";
      (0,-12)*+{x_2}="x2";
      (10,-12)*+{\cdot}="pivot";
      (0,-24)*+{x_3}="x3";
      (10,-24)*+{x_4}="x4";
      (20,-24)*+{x_5}="x5";
      {\ar@{-} "a"; "x1"};
      {\ar@{-} "a"; "x2"};
      {\ar@{-} "a"; "pivot"};
      {\ar@{-} "pivot"; "x3"};
      {\ar@{-} "pivot"; "x4"};
      {\ar@{-} "pivot"; "x5"};
    \end{xy}
    $
    &
    $
    \begin{xy}
      (0,0)*+{};
      (0,-9)*+{=};
    \end{xy}
    $
    &
    $
    \begin{xy}
      (0,0)*+{\cdot}="a";
      (-10,-12)*+{x_1}="x1";
      (0,-12)*+{\cdot}="pivot";
      (10,-12)*+{x_5}="x5";
      (-10,-24)*+{x_2}="x2";
      (0,-24)*+{x_3}="x3";
      (10,-24)*+{x_4}="x4";
      {\ar@{-} "a"; "x1"};
      {\ar@{-} "a"; "pivot"};
      {\ar@{-} "a"; "x5"};
      {\ar@{-} "pivot"; "x2"};
      {\ar@{-} "pivot"; "x3"};
      {\ar@{-} "pivot"; "x4"};
    \end{xy}
    $
    &
    $
    \begin{xy}
      (0,0)*+{};
      (0,-9)*+{=};
    \end{xy}
    $
    &
      $
    \begin{xy}
      (0,0)*+{\cdot}="a";
      (-10,-12)*+{\cdot}="pivot";
      (0,-12)*+{x_4}="x4";
      (10,-12)*+{x_5}="x5";
      (-20,-24)*+{x_1}="x1";
      (-10,-24)*+{x_2}="x2";
      (0,-24)*+{x_3}="x3";
      {\ar@{-} "a"; "x4"};
      {\ar@{-} "a"; "x5"};
      {\ar@{-} "a"; "pivot"};
      {\ar@{-} "pivot"; "x1"};
      {\ar@{-} "pivot"; "x2"};
      {\ar@{-} "pivot"; "x3"};
    \end{xy}
    $
  \end{tabular}
\end{center}

\begin{definition}[Underlying list]
  Let $t \in \Fr_\tensor(\V)$. The \emph{underlying list} of $t$ is
  the word of $\V^*$ defined inductively by
  \[
  U(t) = \begin{cases}
    U(t_1)\of\dots\of U(t_n) & \text{if $t =
      \tensor(t_1,\dots,t_n)$}\\
    t & \text{otherwise}
  \end{cases}
  \]
\end{definition}

\begin{definition}[Left-most bracketing]
  Let $t \in \Fr_\tensor(\V)$. If $U(t) = t_1\of\dots\of t_{n +
    k(n-1)}$, then the left-most bracketing of $t$ is defined
  recursively by
  \[
  \lmb(t_1^{n+k(n-1)}) = \lmb(\tensor(t_1^n),t_{n+1}^{n+ k(n-1)}).
  \]
\end{definition}

We wish to establish that any term  $\Fr_\tensor(\V)$ is equal, in
$\Fr_{C_n}(\V)$, to its left most bracketing. To this end, we define
the rank and the length of a term, which will be useful again for a
related problem in Section \ref{sec:catalancat}.

\begin{definition}\label{def:ranklength}
Let $t \in \Fr_\tensor(\V)$. Define the length of $t$ to be: 
  \[
  L(t) = \begin{cases}
    \sum_{i=1}^n L(t_i) & \text{if $t = \tensor(t_1^n)$}\\
    1 & \text{otherwise}.
  \end{cases}
  \]
  Define the rank, $R(t)$, of $t$ inductively by setting $R(t) = 0$ if
  $t \in \V$ and 
  \[
  R(\tensor(t_1^n)) = \sum_{i=1}^n R(t_i) + \sum_{i=2}^n (i-1)L(t_i) -
  \frac{n(n-1)}{2}.
  \]
  Note that $R(t) = 0$ precisely when $t = \lmb(t)$.
\end{definition}

We may now proceed to show that any term is equivalent to its
left-most bracketing. 

\begin{lemma}\label{lem:lmb}
  For any $t \in \Fr_\tensor(\V)$, we have $t =_{C_n} \lmb(t)$
\end{lemma}
\begin{proof}
 Let $t \in \Fr_\tensor(\V)$. Let $R(t)$ and $L(t)$ be as in
 Definition \ref{def:ranklength}. We proceed by double induction on
 $R(t)$ and $L(t)$ to show that $t=_{C_n} \lmb(t)$. If $L(t) = 1$ then
 the statement is trivial. We also have that $R(t) = 0$ if and only if
 $t = \lmb(t)$. 

 Suppose that $L(t) > 1$ and $R(t) > 0$, so that $t =
 \tensor(t_1^n)$. Let $i$ be the greatest integer with the property that
 $t_i \notin \V$. If $i = 1$, then $t = \lmb(t)$ by induction on
 $L(t)$. If $i > 1$, then $t_i = \tensor(u_1^n)$ and set
 \[
 t' = \tensor(t_1^{i-2}, \tensor(t_{i-1}, u_1^{n-1}),u_n,t_{i+1}^n).
 \]
 A single application of one of the equations in $C_n$ establishes
 that $t =_{C_n} t'$. Since $R(t) - R(t') = \sum_{i=1}^{n-1} L(u_i)$,
 we have $R(t') < R(t)$ and the statement follows by induction on
 $R(t)$.
\end{proof}

In order to manipulate elements of $\Struct(C_n)$ effectively, we
introduce the notion of a seed.

\begin{definition}[Seed]
  Let $\F$ be a graded set of function symbols on some set $\V$ and let
  $\rho$ be a partial function $\Fr_\F(\V) \to \Fr_\F(\V)$. A
  \emph{seed} for $\rho$ is a pair of terms $s,t \in \Fr_\F(\V)$ such
  that the graph of $\rho$ is equal to
  $\{(s^\varphi,t^\varphi)~|~\varphi \in [\V,\Fr_\F(\V)]\}$.
\end{definition}

In particularly nice cases, we can construct seeds for any operator in
a structure monoid. 

\begin{lemma}[Dehornoy \cite{Dehornoy:braids}]\label{lem:seeds}
  Let $\T$ be a balanced equational theory that contains precisely one
  function symbol. Then, each operator $\rho \in \Struct(\T)$ admits a
  seed. 
\end{lemma}

It follows from Lemma \ref{lem:expansion} that $C_n$ is composable
and we may, therefore, form the group $\Struct_G(C_n)$. In order to
facilitate the passage from members of $\Struct_G(C_n)$, to members of
$F_{n,1}$, we introduce the tree generated by a term. 

\begin{definition}
  For a term $t \in \Fr_\tensor(\V)$, let $T(t)$ denote the $n$-ary
  tree obtained via the following construction:
  
  \begin{itemize}
    \item If $t = \tensor(t_1,\dots,t_n)$, then $T(t)$ is equal to:
      \[
      \begin{xy}
        (0,0)*+{\cdot}="a";
        (-18,-14)*+{T(t_1)}="t1";
        (-7,-14)*+{T(t_2)}="t2";
        (2,-14)*+{\dots}="dots";
        (12,-14)*+{T(t_n)}="tn";
        {\ar@{-} "a"; "t1"};
        {\ar@{-} "a"; "t2"};
        {\ar@{-} "a"; "tn"};
      \end{xy}
      \]
    \item Otherwise, $T(t) = \cdot$
\end{itemize}

\end{definition}

\begin{theorem}\label{thm:Fn1}
  $\Struct_G(C_n) \cong F_{n,1}$.
\end{theorem}
\begin{proof}
We denote the seed of $\rho \in \Struct_G(C_n)$, which exists by Lemma
\ref{lem:seeds}, by $(s_\rho,t_\rho)$. We claim that the following map
is an isomorphism:
\begin{eqnarray*}
    \Struct_G(C_n) &\stackrel{\theta}{\longrightarrow}& F_{n,1}\\
    \rho &\longmapsto& [(T(s_\rho),T(t_\rho),id)]
\end{eqnarray*}

It is routine to see that $\Theta$ is a homomorphism. Suppose that
$\rho,\rho' \in \Struct_G(C_n)$ and that $\Theta(\rho) =
\Theta(\rho')$. It follows that $\rho$ and $\rho'$ have the same seed,
so $\rho = \rho'$ and $\Theta$ is faithful.

By Lemma \ref{lem:capture}, in order to establish
that $\Theta$ is surjective, we need only show that $t_1 =_{C_n} t_2$
whenever $t_1,t_2 \in \Fr_\tensor(\V)$ and $U(t_1) = U(t_2)$. By
Lemma \ref{lem:lmb}, we have $t_1 =_{C_n} \lmb(t_1) =_{C_n} \lmb(t_2)
=_{C_n} t_2$, so $\Theta$ is surjective and, hence, an isomorphism.
\end{proof}

\subsection{Symmetric Catalan Algebras and $G_{n,1}$} We saw in
Section \ref{sec:groups} that the leaves of a tree may be ordered by
the lexicographic ordering on their addresses. An $n$-ary tree symbol
$[(T_1,T_2,\sigma)]$ may thereby be viewed as a pair of tree diagrams,
together with a permutation of the leaves of $T_1$. Thus, in order to
obtain an equational theory whose structure group is $G_{n,1}$ we need
to add the ability to arbitrary permute variables to Catalan
algebras. Recalling that the symmetric group is generated by
transpositions of adjacent elements, we are led to the following
definition.

\begin{definition}[Symmetric $n$-Catalan Algberas]
  The theory of symmetric $n$-catalan algebras extends that of
  $n$-catalan algebras with the following equations, where $1 \le i <
  n$:
  \[
  \tensor(x_1^{i-1},x_i,x_{i+1},x_{i+2}^n) =
  \tensor(x_1^{i-1},x_{i+1},x_i,x_{i+2}^n).         
  \]
  We denote the theory of symmetric $n$-catalan algebras by $SC_n$.
\end{definition}

Symmetric $n$-catalan algebras essentially add an action of the
symmetric group on the indices of $\tensor$. In general, this is
sufficient to induce an action of a symmetric group on the variables
of any term in $\Fr_\tensor(\V)$. In the binary case, we
recover the definition of commutative semigroups.

\begin{theorem}\label{thm:Gn1}
  $\Struct_G(SC_n) \cong G_{n,1}$. 
\end{theorem}
\begin{proof}
  For $\rho \in \Struct_G(SC_n)$, let $(s_\rho,t_\rho)$ represent its
  seed, which exists by Lemma \ref{lem:seeds}. Since $SC_n$ is linear,
  $s_\rho$ and $t_\rho$ are linear and $\supp(s_\rho) =
  \supp(t_\rho)$. Let $\pi(\rho)$ be the permutation of
  $\supp(s_\rho)$ induced by the permutation $U(s_\rho) \to
  U(t_\rho)$. Consider the following map:
  \begin{eqnarray*}
    \Struct_G(C_n) &\stackrel{\theta}{\longrightarrow}& G_{n,1}\\
    \rho &\longmapsto& [(T(s_\rho),T(t_\rho),\pi(\rho))]
  \end{eqnarray*}
   A similar argument to the proof of Theorem \ref{thm:Fn1}
  establishes that $\Theta$ is an isomorphism.
\end{proof}

We now know that $F_{n,1}$ and $G_{n,1}$ are the structure groups of
catalan algebras and of symmetric catalan algebras, respectively. We also
know that if we can construct coherent categorifications of these
algebras, then we can apply Theorem \ref{thm:group} to obtain
presentations of these groups. In the following section, we set about
the task of constructing a coherent categorification of catalan
algebras.

\section{Catalan categories and $F_{n,1}$}\label{sec:catalancat} 

In order to obtain a presentation for $\Struct_G(C_n)$ and, hence, for
$F_{n,1}$ along the lines of that provided by Dehornoy for $F$
\cite{Dehornoy:thompson}, we need to obtain a coherent
categorification of $C_n$. The immediate problem is discerning a set
of diagrams whose commutativity imply the commutativity of all
diagrams in $\Fr_{\widehat{C_n}}(\hV)$. As we shall see in this
section, the following definition suffices for this purpose. While the
coherence axioms that we have chosen may seem slightly cryptic, the
reason for their choice will become apparent in the proof that the
resulting categorification is coherent. We shall make frequent use of
the following useful shorthand: For $1 \le i \le n$ and a morphism
$\rho:t_i \to t_i'$, we set   
  \[
  \tensor^i(\rho) = \tensor(1_{t_1},\dots,
  1_{t_{i-1}},\rho,1_{t_{i+1}},\dots,1_{t_n}).
  \]

\begin{definition}
  A discrete $n$-catalan category is the categorification of
  $C_n$ that consists of:
  \begin{itemize}
  \item A discrete category $\hV$.
  \item A functor $\tensor: \hV^n \to \hV$.
  \item For $1 \le i < n$, a natural
    isomorphism $\alpha_i$ with the following components:
    \[
    \alpha_i(x_1^{2n-1}) :
    \tensor(x_1^i,\tensor(x_{i+1}^{i+n}),x_{i+n+1}^{2n-1}) 
    \to
    \tensor(x_1^{i-1},\tensor(x_i^{i+ n-1}), x_{i+n}^{2n-1})
    \]
  \end{itemize}

  \noindent \textbf{Pentagon axiom:} For $1 \le i \le n-1$, the following 
    diagram commutes, where $X = x_1^{i-1}$ and  $Z = z_1^{n-i-1}$:  
    \[
    \begin{xy}
      (0,4)*+{{\tensor(X,y_1,\tensor(y_2^n,\tensor(y_{n+1}^{2n})),Z)}}="1";
      (-33,-15)*+{\tensor(X,\tensor(y_1^n),\tensor(y_{n+1}^{2n}),Z)}="2";
      (30,-15)*+{\tensor(X,y_1,\tensor(y_i^{n-1},\tensor(y_n^{2n-1}),y_{2n}),Z)}
      ="3"; 
      (-33,-35)*+{\tensor(X,\tensor(\tensor(y_1^n),y_{n+1}^{2n-1}),y_{2n},Z)}
      ="4";
      (33,-35)*+{\tensor(X,\tensor(y_1^{n-1},\tensor(y_n^{2n-1})),y_{2n},Z)}
      ="5";
      {\ar@{->}_<<<<<<<{\alpha_i}@/_/ "1"; "2"}
      {\ar@{->}^<<<<<<<{\tensor^{i+1}(\alpha_{n-1})}@/^/ "1"; "3"}
      {\ar@{->}^{\alpha_i} (26,-19); (26,-32)}
      {\ar@{->}_{\alpha_i} "2"; "4"}
      {\ar@{->}^{\tensor^i(\alpha_{n-1}\of\dots\of\alpha_1)} @/^1.8pc/
        (22,-38); (-33,-39)} 
    \end{xy}
    \]
    \medskip\\
    \bigbreak
    \noindent\textbf{Adjacent associativity axiom:} For $1 \le i \le
    n-2$ , the following 
    diagram commutes, where $X = x_1^{i-1}$ and  $Z = z_1^{n-i-2}$: 
    
    \[
    \begin{xy}
      (10,4)*+{\tensor(X,y_1,\tensor(y_2^{n+1}),\tensor(y_{n+2}^{2n+1}),Z)}="a";
      (-28,-12)*+{\tensor(X,\tensor(y_1^n),y_{n+1},\tensor(y_{n+2}^{2n+1}),
        Z)}="b"; 
      (28,-24)*+{\tensor(X,y_1,\tensor(\tensor(y_2^{n+1}),y_{n+2}^{2n}),
        y_{2n+1},Z)}="c";
      (-32,-35)*+{\tensor(X,\tensor(y_1^n),\tensor(y_{n+1}^{2n-1}),y_{2n}^{2n+1},
      Z)}="d";
      (28,-46)*+{\tensor(X,\tensor(y_1,\tensor(y_2^{n+1}),y_{n+2}^{2n-1}),
        y_{2n}^{2n-1},Z)}="e";
      (-22,-59)*+{\tensor(X,\tensor(\tensor(y_1^n),y_{n+1}^{2n-1}),y_{2n}^{2n+1},
      Z)}="f";
    {\ar@{->}_<<<<<<<{\alpha_i}@/_1pc/ "a"; (-29,-11)}
    {\ar@{->}^{\alpha_{i+1}}@/^/ "a"; (28,-20)}
    {\ar@{->}_{\alpha_{i+1}} (-32,-18); "d"}
    {\ar@{->}^{\alpha_i} (27,-28); (27,-41)}
    {\ar@{->}_{\alpha_i}  (-35,-38);(-35,-55)}
    {\ar@{->}^{\tensor^i(\alpha_1)}@/^/ (23,-49); (1,-59)}
  \end{xy}
  \]

  We denote the theory of discrete $n$-catalan categories by $\Cn$ and
  the free $\Cn$ category on $\hV$ by $\Fr_\Cn(\hV)$.
\end{definition} 

In the case where $n=2$, the pentagon axiom reduces to Mac Lane's
pentagon axiom for monoidal categories from Example \ref{ex:monoidal}
and the adjacent associativity axiom is empty, so we recover the usual
definition of a coherently associative bifunctor. 

\begin{definition}[Positive/Negative]
  A singular morphism of $\Fr_\Cn(\hV)$ that contains an instance of $\alpha_i$
  is called \emph{positive} and one that contains an instance of
  $\alpha_i^{-1}$ is called \emph{negative}. A morphism in
  $\Fr_\Cn(\hV)$ is called \emph{positive} if it is an identity or a
  composite of positive morphisms and \emph{negative} if it is a
  composite of negative morphisms. 
\end{definition}

It follows from the proof of Lemma \ref{lem:lmb} that there is always a
positive morphism $t \to \lmb(t)$ in $\Fr_\Cn(\hV)$. 
In order to show that $\Cn$ is a
coherent categorification of $C_n$, we need to show that any diagram
built out of the singular and identity morphisms of $\Cn$ commutes. As
our first step towards this goal, we show that there is a unique
positive morphism $t \to \lmb(t)$. 

\begin{lemma}\label{lem:unique}
  Let $t \in \Ob(\Fr_\Cn(\hV))$. There is a unique positive morphism
  $t \to \lmb(t)$ in $\Fr_\Cn(\hV)$. 
\end{lemma}
\begin{proof}
  Let $t \in \Ob(\Fr_\Cn(\hV))$. It follows from the proof of Lemma
  \ref{lem:lmb} that there is a positive morphism $t \to
  \lmb(t)$. Suppose that $\varphi,\psi:t \to \lmb(t)$,
  that $\varphi = \varphi_1 \of \varphi_2$ and that $\psi = \psi_1 \of
  \psi_2$. By Lemma \ref{lem:singular}, we may assume that $\varphi_1$
  and $\psi_1$ are singular.  Let $R(t)$ 
  and $L(t)$ be defined as in Definition \ref{def:ranklength}. We
  proceed by double induction on $R(t)$ and $L(t)$ to show that there
  exists an object $w$ in $\Fr_\Cn(\hV)$ making the following diagram
  commute:
  \begin{displaymath}
  \vcenter{
    \def\objectstyle{\scriptstyle}
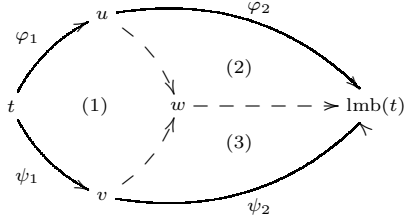
    \begin{xy}
      (0,0)*+{t}="s";
      (12,12)*+{u}="u1";
      (12,-12)*+{v}="u2";
      (22,0)*+{w}="t";
      (48,0)*+{\lmb(t)}="ns";
      (11,0)*+{(1)};
      (30,5)*+{(2)};
      (30,-5)*+{(3)};
      {\ar@{->}@/^0.4pc/^{\varphi_1} "s"; "u1"};
      {\ar@{->}@/_0.4pc/_{\psi_1} "s"; "u2"};
      {\ar@{-->}@/^0.4pc/ "u1"; "t"};
      {\ar@{-->}@/_0.4pc/ "u2"; "t"};
      {\ar@{->}@/^1.3pc/^{\varphi_2} "u1"; "ns"};
      {\ar@{->}@/_1.3pc/_{\psi_2} "u2"; "ns"};
      {\ar@{-->} "t"; "ns"};
    \end{xy}
    }
  \end{displaymath}
  By induction on $R(t)$, we have that the subdiagrams labelled $(2)$ and
  $(3)$ commute. So, we need only establish the existence and
  commutativity of the subdiagram labelled $(1)$. If $R(t) = 0$ or
  $L(t) = 1$, then the statement is trivial, so suppose that $R(t) >
  0$ and $L(t) > 1$. If $\varphi_1 = \psi_1$, then take $w = u = v$
  and subdiagram $(1)$ commutes trivially. 

  Suppose that $\varphi_1 \ne \psi_1$. Then, either $\varphi_1 =
  \tensor^i(\varphi_1')$ or $\varphi_1 = \alpha_i(t_1^n)$ and there
  are similar possibilities for $\psi_1$. We proceed by case analysis
  on the form of $\varphi_1$ and $\psi_1$.
   
   Suppose that $\varphi_1 = \tensor^i(\varphi_1')$ and $\psi_1 =
   \tensor^j(\psi_1')$. If $i = j$, then the whole diagram commutes by
   induction on $L(t)$. Suppose that $i \ne j$. Then, without loss of
   generality, $i < j$ and we may take $(1)$ to be the following
   diagram, which commutes by the functoriality of $\tensor$:
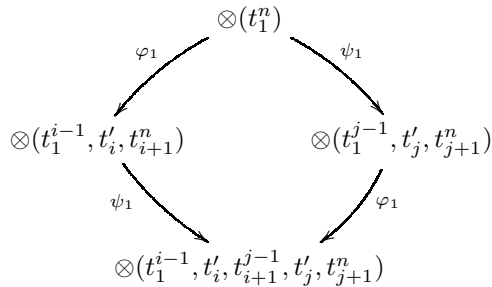
   \[
   \begin{xy}
     (0,0)*+{\tensor(t_1^n)}="1";
     (-20,-16)*+{\tensor(t_1^{i-1},t_i',t_{i+1}^n)}="2";
     (20,-16)*+{\tensor(t_1^{j-1},t_j',t_{j+1}^n)}="3";
     (0,-33)*+{\tensor(t_1^{i-1},t_i',t_{i+1}^{j-1},
       t_j',t_{j+1}^n)}="4";
     {\ar@{->}_{\varphi_1}@/_/ "1"; "2"}
     {\ar@{->}^{\psi_1}@/^/ "1"; "3"}
     {\ar@{->}_<<<<{\psi_1}@/_/ "2"; "4"}
     {\ar@{->}^<<<<{\varphi_1}@/^/ (15,-20); "4"}
   \end{xy}
   \]

   Suppose that $\varphi_1 = \tensor^i(\varphi_1')$ and $\psi_1 =
   \alpha_j(t_1^n)$. If $i \ne j$, then without loss of generality $i
   < j$ and $t_j = \tensor(u_1^n)$. If $i \ne j-1$, then we may take
   $(1)$ to be the following square, which commutes by the naturality
   of $\alpha_i$: 
      \[
   \begin{xy}
     (0,0)*+{\tensor(t_1^{j-1},\tensor(u_1^n),t_{j+1}^n)}="1";
     (-30,-16)*+{\tensor(t_1^{i-1},t_i',t_{i+1}^{j-1},
       \tensor(u_1^n),t_{j+1}^n)}="2";  
     (30,-16)*+{\tensor(t_1^{j-2},\tensor(t_{j-1},u_1^{n-1}),u_n,
       t_{j+1}^n)}="3";  
     (0,-33)*+{\tensor(t_1^{i-1},t_i',t_{i+1}^{j-2},
       \tensor(t_{j-1},u_1^{n-1}),u_n,t_{j+1}^n)}="4"; 
     {\ar@{->}_{\varphi_1}@/_/ "1"; "2"}
     {\ar@{->}^{\psi_1}@/^/ "1"; "3"}
     {\ar@{->}_<<<<{\psi_1'}@/_/ "2"; "4"}
     {\ar@{->}^<<<<{\varphi_1'}@/^/ (23,-20); "4"}
   \end{xy}
   \]
   
   If $i = j-1$, then we may take $(1)$ to be a similar naturality
   square. Suppose that $i=j$. If $\varphi_1' =
   \tensor^k(\varphi_1'')$, then we may take $(1)$ to be a naturality
   square. If $\varphi_1' = \alpha_k(u_1^n)$ then we have two
   cases. If $k \ne (n-1)$, then we may take $(1)$ to be a naturality
   square. If $k = n-1$, then we may take $(1)$ to be an instance of
   the pentagon axiom, which commutes by assumption.

   Finally, we are left with the case where $\varphi_1 =
   \alpha_i(t_1^n)$ and $\psi_1 = \alpha_j(t_1^n)$. If $|i-j| > 1$,
   then we may take $(1)$ to be a naturality diagram. If $|i-j| = 1$,
   then we may take $(1)$ to be an instance of the adjacent
   associativity axiom. 
\end{proof}

We now know that every object in $\Fr_\Cn(\hV)$ has a unique positive
morphism to its left-most bracketing. With a little work, we can
bootstrap this result in order to show that there is a unique
morphism - positive, negative or otherwise - between any two arbitrary
objects in $\Fr_\Cn(\hV)$.

\begin{theorem}\label{thm:coherencecatalan}
  $\Cn$ is a coherent categorification of $C_n$.
\end{theorem}
\begin{proof}
 Suppose that $\varphi:s \to t$ is a reduction in $\Fr_\Cn(\hV)$. By
 Lemma \ref{lem:singular}, 
\[
\varphi = s \stackrel{\varphi_1}{\to} s_1 \stackrel{\varphi_2}{\to}
s_2 \to \dots \stackrel{\varphi_{n-1}}{\to} s_{n-1}
\stackrel{\varphi_n}{\to} t, 
\]
where each $\varphi_i$ is singular. By Lemma \ref{lem:unique}, each
term  $t \in \Fr_\Cn(\hV)$ has a unique map $N_t:t \to \lmb(t)$. We
claim that each rectangle in the following diagram commutes: 
\[
\vcenter{
  \xymatrix@=0.99pc{
    {s} \ar[r]^{\varphi_1} \ar[dd]_{N_s} & {s_1} \ar[r]^{\varphi_2}
    \ar[dd]_{N_{s_1}} & {s_2} \ar[r]^{\varphi_3} \ar[dd]_{N_{s_2}} &
    {\dots} \ar[r]^{\varphi_{n-1}} & {s_{n-1}} \ar[r]^{\varphi_n}
    \ar[dd]_{N_{s_{n-1}}} & {t} \ar[dd]_{N_t}\\
    \\
    {\lmb(s)} \ar@{=}[r] & {\lmb(s_1)} \ar@{=}[r] & {\lmb(s_2)} \ar@{=}[r] &
    {\dots} \ar@{=}[r] & {\lmb(s_{n-1})} \ar@{=}[r] & {\lmb(t)} 
  }
}
\]
If $\varphi_i$ is positive, then it follows immediately from Lemma
\ref{lem:unique} that $\varphi_i\of N_{s_i} = N_{s_{i-1}}$. If
$\varphi_i$ is negative, then Lemma \ref{lem:unique} implies that
$\varphi_i^{-1}\of N_{s_{i-1}} = N_{s_i}$, which implies that
$\varphi_i\of N_{s_i} = N_{s_{i-1}}$. Since each rectangle commutes,
we have $\varphi\of N_t = N_s$, which implies that $\varphi = N_s\of
N_t^{-1}$. Since $N_s$ and $N_t$ are unique and we did not rely on a
particular choice of $\varphi$, we conclude that $\Cn$ is coherent.
\end{proof}

With Theorem \ref{thm:coherencecatalan} in hand, we can obtain a
presentation for $F_{n,1}$, which generalises the presentation for $F$
given in \cite{Dehornoy:thompson}. 

\begin{corollary}
  $\P_G(\Cn) \cong F_{n,1}$
\end{corollary}
\begin{proof}
  By Theorem \ref{thm:coherencecatalan} and Theorem \ref{thm:group},
  we have $\P-G(\Cn) \cong \Struct_G(C_n)$. It follows then from Theorem
  \ref{thm:Fn1} that $\P_G(\Cn) \cong F_{n,1}$.
\end{proof}

In the following section, we shall obtain a coherent categorification
of $SC_n$ and, thereby, a presentation of $G_{n,1}$.

\section{Symmetric catalan categories and
  $G_{n,1}$}\label{sec:scatalancat} 

Our goal in this section is to construct a coherent categorification
of symmetric catalan algebras. The coherence theorem for catalan
categories, Theorem \ref{thm:coherencecatalan}, reduces this problem
to ensuring that any two sequences of transpositions of the objects
appearing in a term realise the same permutation. In other words, our
categorification needs to somehow encode a presentation of the
symmetric group whose generators correspond to transpositions of
adjacent variables. Such a presentation is well known, having been
constructed by Moore \cite{Moore:presentation}. This presentation has
generators $T_1,\dots, T_{n-1}$ and the following relations:
\begin{eqnarray*}
  T_i^2 = 1 & \text{for $1 \le i \le n-1$}\\
  (T_iT_{i+1})^3 = 1 & \text{for $1 \le i \le n-2$}\\
  (T_iT_k)^2 = 1 &\text{for $1 \le i \le k-2$}
\end{eqnarray*}
With this presentation in mind, we may now construct a reasonable
categorification of $SC_n$. Recall our shorthand that for $1 \le i \le
n$ and a morphism $\rho:t_i \to t_i'$, we have
  \[
  \tensor^i(\rho) = \tensor(1_{t_1},\dots,
  1_{t_{i-1}},\rho,1_{t_{i+1}},\dots,1_{t_n}).
  \]

\begin{definition}
  For $n \ge 2$, a discrete symmetric $n$-catalan category, is
  a discrete $n$-catalan category on the category $\hV$, together
  with, for $1 \le i \le n-1$, a natural isomorphism $\tau_i$ with
  components 
  \[
  \tau_i(t_1^n) : \tensor(t_1^{i-1},t_i,t_{i+1},t_{i+2}^n) \to
  \tensor(t_1^{i-1},t_{i+1},t_{i},t_{i+2}^n),
  \]
  satisfying the following axioms:\\
  
  \noindent\textbf{Involution axiom:} For $1 \le i \le n-1$, the
  following diagram commutes:
  \[
  \begin{xy}
    (0,0)*+{\tensor(t_1^n)}="1";
    (0,-20)*+{\tensor(t_1^n)}="2";
    (45,-20)*+{\tensor(t_1^{i-1},t_{i+1},t_i,t_{i+2}^n)}="3";
    {\ar@{->}_{1} "1"; "2"}
    {\ar@{->}^{\tau_i} "1"; "3"}
    {\ar@{->}^>>>>>>>>{\tau_i} "3"; "2"}
  \end{xy}
  \]
   \medskip
  \noindent\textbf{Compatibility axiom:} For $2 \le i \le n$ and $1
  \le j \le n-2$, the following
  diagram commutes, where $W = w_1^i$ and $Z = z_1^{n-i}$: 
  \[
  \begin{xy}
    (0,0)*+{\tensor(W,x,\tensor(y_1^n),Z)}="1";
    (-30,-20)*+{\tensor(W,\tensor(x,y_1^{n-1}),y_n,Z)}="2";
    (33,-20)*+{\tensor(W,x,\tensor(y_1^{j-1},y_{j+1},y_j,y_{j+2}^n),Z)}
    ="3";
    (0,-42)*+{\tensor(W,\tensor(x,y_1^{j-1},y_{j+1},y_j,y_{j+2}^{n-1}),
      y_n,Z)}="4";
    {\ar@{->}_{\alpha_{i-1}}@/_/ "1"; "2"}
    {\ar@{->}^{\tensor^i(\tau_j)}@/^/ "1"; "3"}
    {\ar@{->}_<<<<<<{\tensor^{i-1}(\tau_{j+1})}@/_/ "2"; "4"}
    {\ar@{->}^<<<<<<{\alpha_{i-1}}@/^/ (27,-23); "4"}
  \end{xy}
  \]

  \noindent\textbf{$3$-cycle axiom:} For $1 \le i \le n-2$, the
  following diagram commutes:
  \[
  \begin{xy}
    (0,0)*+{\tensor(t_1^n)}="1";
    (-30,-20)*+{\tensor(t_1^{i-1},t_{i+1},t_i,t_{i+2}^n)}="2";
    (30,-20)*+{\tensor(t_1^i,t_{i+2},t_{i+1},t_{i+3}^n)}="3";
    (-30,-40)*+{\tensor(t_1^{i-1},t_{i+1},t_{i+2},t_{i},t_{i+3}^n)}="4";
    (30,-40)*+{\tensor(t_1^{i-1},t_{i+2},t_{i},t_{i+1},t_{i+3}^n)}="5";
    (0,-60)*+{\tensor(t_1^{i-1},t_{i+2},t_{i+1},t_{i},t_{i+3}^n)}="6";
    {\ar@{->}_{\tau_i}@/_/ "1"; "2"}
    {\ar@{->}^{\tau_{i+1}}@/^/ "1"; "3"}
    {\ar@{->}_{\tau_{i+1}} "2"; "4"}
    {\ar@{->}^{\tau_i} (26,-24); (26,-37)}
    {\ar@{->}_<<<<<<<{\tau_{i}}@/_/ (-32,-44); (-15,-56)}
    {\ar@{->}^<<<<<<<{\tau_{i+1}}@/^/ (22,-45); "6"}
  \end{xy}
  \]

  \noindent\textbf{Hexagon axiom:} For $1 \le i \le n-1$, the
  following diagram commutes, where $W = w_1^{i-1}$ and $Z =
  z_1^{n-i-1}$:
  \[
  \begin{xy}
    (0,0)*+{\tensor(W,\tensor(x_1^n),y,Z)}="1";
    (-30,-20)*+{\tensor(W,y,\tensor(x_1^n),Z)}="2";
    (30,-20)*+{\tensor(W,x_1,\tensor(x_2^n,y),Z)}="3";
    (-30,-40)*+{\tensor(W,\tensor(y,x_1^{n-1}),x_n,Z)}="4";
    (30,-40)*+{\tensor(W,x_1,\tensor(y,x_2^n),Z)}="5";
    (0,-60)*+{\tensor(W,\tensor(x_1,y,x_2^{n-1}),x_n,Z)}="6";
    {\ar@{->}_{\tau_i}@/_/ "1"; "2"}
    {\ar@{->}^{\alpha_i^{-1}}@/^/ "1"; "3"}
    {\ar@{->}_{\alpha_i} "2"; "4"}
    {\ar@{->}^{\tensor^{i+1}(\tau_{n-1}\of\dots\of\tau_1)} (26,-24);
      (26,-37)} 
    {\ar@{->}_<<<<<<<{\tensor^i(\tau_1)}@/_/ (-32,-44); (-15,-56)}
    {\ar@{->}^<<<<<<<{\alpha_i}@/^/ (22,-45); "6"}
  \end{xy}
  \]
  We denote the theory of discrete symmetric $n$-catalan categories
  by $\SCn$ and the free $\SCn$-category on $\hV$ by $\Fr_\SCn(\hV)$.
\end{definition}

The hexagon axiom ensures that we may replace a transposition of the
form $\tau_i(t_1^{i-1},\tensor(u_1^n),t_{i}^{n-1})$ with a sequence of
transpositions involving only the terms $t_1^{n-1}$ and
$u_1^n$. One might posit the commutativity of a diagram that serves the
same purpose for a morphism of the form
$\tau_i(t_1^{i},\tensor(u_1^n),t_{i+1}^{n-1})$. Doing so leads to the
\emph{dual hexagon diagram}, which has the following form, where $2
\le i \le n$ and $W = w_1^{i-2}$ and $Z = z_1^{n-i}$: 
\[
\begin{xy}
  (0,0)*+{\tensor(W,x,\tensor(y_1^n),Z)}="1";
  (-30,-20)*+{\tensor(W,\tensor(y_1^n),x,Z)}="2";
  (30,-20)*+{\tensor(W,\tensor(x,y_1^{n-1}),y_n,Z)}="3";
  (-30,-40)*+{\tensor(W,y_1,\tensor(y_2^n,x),Z)}="4";
  (30,-40)*+{\tensor(W,\tensor(y_1^{n-1},x),y_n,Z)}="5";
  (0,-60)*+{\tensor(W,y_1,\tensor(y_2^{n-1},x,y_n),Z)}="6";
    {\ar@{->}_{\tau_i}@/_/ "1"; "2"}
    {\ar@{->}^{\alpha_i}@/^/ "1"; "3"}
    {\ar@{->}_{\alpha_i^{-1}} "2"; "4"}
    {\ar@{->}^{\tensor^{i}(\tau_{1}\of\dots\of\tau_{n-1})} (26,-24);
      (26,-37)} 
    {\ar@{->}_<<<<<<<{\tensor^{i+1}(\tau_{n-1})}@/_/ (-32,-44); (-15,-56)}
    {\ar@{->}^<<<<<<<{\alpha_i^{-1}}@/^/ (22,-45); "6"}
\end{xy}
\]

\begin{lemma}\label{lem:dualhex}
  The dual hexagon diagram commutes in $\Fr_\SCn(\hV)$.
\end{lemma}
\begin{proof}
  By the involution axiom, we have
  \[\tau_i(t_1^{i},\tensor(u_1^n),t_{i+1}^{n-1}) =
  [\tau_i(t_1^{i-1},\tensor(u_1^n),t_{i}^{n-1})]^{-1}.\] Using the
  hexagon axiom and ignoring component labels, we have:
  \begin{eqnarray*}
    \tau_i &=& \tau_i^{-1}\\
    &=& (\alpha_i^{-1}\of\tensor^{i+1}(\tau_{n-1}\of\dots\of\tau_1)
    \of \alpha_i \of\tensor^i(\tau_1)\of\alpha_i^{-1})^{-1}
  \end{eqnarray*}
  So, we have:
  \begin{eqnarray}
    \nonumber \tau_i\of\alpha_i^{-1}\of\tensor^{i+1}(\tau_{n-1}) &=&
    \alpha_i\of\tensor^{i}(\tau_1)\of\alpha_i^{-1}\of
  \tensor^{i+1}(\tau_1,\dots,\tau_{n-2})\\  
    &=& \alpha_i\of\tensor^{i}(\tau_1)\of\alpha_i^{-1}\of
    \tensor^{i+1}(\tau_1)\of\dots\of\tensor^{i+1}(\tau_{n-2})\label{eqn:lhs}
  \end{eqnarray}
  From the compatibility axiom, we have $\tensor^{i+1}(\tau_j) =
  \alpha_i\of\tensor^i(\tau_{j+1})\of\alpha_i^{-1}$. This implies:
  \begin{eqnarray*}
    (\ref{eqn:lhs}) &=& \alpha_i\of\tensor^i(\tau_1)\of\dots
    \of\tensor^i(\tau_{n-1})\of\alpha_i^{-1}\\
    &=& \alpha_i\of\tensor^i(\tau_1\of\dots\of\tau_{n-1})
    \of\alpha_i^{-1}
  \end{eqnarray*}
  Therefore, the dual hexagon diagram commutes in $\Fr_\SCn(\hV)$.
\end{proof}

In the $n=2$ case, the axiomatisation of $\SCn$ reduces to the theory
of a coherently associative and commutative bifunctor given in Example
\ref{ex:monoidal}. The main
result of this section establishes that $\SCn$ is a suitable
generalisation of this case. 

\begin{theorem}\label{thm:coherencescatalan}
  $\SCn$ is a coherent categorification of $SC_n$. 
\end{theorem}
\begin{proof}
  By Theorem \ref{thm:coherencecatalan}, we may assume that all of the
  associativity maps are strict equalities. Thus, an object of
  $\Fr_\SCn(\hV)$ may be represented as $\tensor(t_1^m)$, where each
  $t_i$ is an object in $\hV$ and $m = n + k(n-1)$, for some $k \ge
  0$. Lemma   \ref{lem:dualhex} and the hexagon axiom imply that it
  suffices to   consider transpositions of adjacent variables that are
  objects of   $\hV$. So, for a given object $t := \tensor(t_1^m)$, we
  need only   consider the $m-1$ induced transposition natural
  isomorphisms  
  \[
    T_i(t_1^m) : \tensor(t_1^{i-1},t_i,t_{i+1},t_{i+2}^m) \to
  \tensor(t_1^{i-1},t_{i+1},t_{i},t_{i+2}^m).
  \]
  In order to establish coherence, we have to show that every
  permutation of $t_1^m$ is unique. That is, we have to show that the
  induced transposition maps satisfy the defining relations for the
  symmetric group of order $m$. 

  The compatibility axiom implies that each $T_i$ is unique.  By the
  naturality of the maps $T_i$, we have $T_i\of T_k = 
  T_k\of T_i$ for all $1 \le i \le k-2$. The involution axiom implies
  that $T_i^2 = 1$. Thus, it only remains to establish that $(T_i\of
  T_{i+1})^3 = 1$. For $n=2$, we may use the proof from Mac Lane
  \cite{MacLane_natural}. Suppose that $n \ge 3$. Since the
  associativity maps are taken to be strict equalities, we may assume
  that $t$ has  the form
  $\tensor(R,\tensor(S,t_i,t_{i+1},t_{i+2},U),V)$, where $R,S,U$ and
  $V$ are sequences of objects of $\hV$. The result then follows from
  the $3$-cycle axiom. 
\end{proof}

With the coherence theorem in hand, we can construct a presentation of
\newline $\Struct_G(SC_n)$ and, therefore, of $G_{n,1}$, which generalises the
presentation for $V$ given in \cite{Dehornoy:thompson}.

\begin{corollary}
  $\P_G(\SCn) \cong G_{n,1}$
\end{corollary}
\begin{proof}
  By Theorem \ref{thm:coherencescatalan} and Theorem \ref{thm:group},
  we have \[\P_G(\SCn) \cong \Struct_G(SC_n).\] It follows then
  from Theorem 
  \ref{thm:Gn1} that $\P_G(\Cn) \cong G_{n,1}$.
\end{proof}

\section{Conclusions and further work}\label{sec:conclusion}

We have demonstrated, by way of Theorem \ref{thm:monoid} and Theorem
\ref{thm:group}, that there is a close relationship between structure
monoids and coherent categorical theories. This relationship is quite
powerful, as illustrated by the fact that we were able to exploit it
in order to obtain new presentations of $F_{n,1}$ and $G_{n,1}$. 

While we only dealt with
invertible categorical structures, it is straightforward to extend the
constructions of Section \ref{sec:free} to structures involving a mix
of invertible and non-invertible natural transformations. Within this
setting, it is possible to develop an abstract coherence theorem that
applies to a large array of structures and, inter alia, yields
presentations of a wide variety of structure monoids. A general
coherence theorem along the lines of the proof of Theorem
\ref{thm:coherencecatalan} is developed in \cite{CohenJohnson:nf}. A
more powerful, though more difficult to apply, general coherence
theorem applying mainly to badly behaved non-invertible structure is
developed in \cite{Cohen:nnf}. 

The presentations that arise from coherence theorems are all
infinite, regardless of whether or not the associated algebraic
structure is finitely presentable. In particular, this is the case for
the presentations of $F_{n,1}$ and $G_{n,1}$ that arise from the
coherence theorems for catalan categories and symmetric catalan
categories. However, these presentations arise from finitely presented
categorical structures and so retain some amount of finiteness.  Conditions on
a finitely presented categorical structure ensuring that it is
coherent via finitely many coherence conditions are obtained in
\cite{Cohen:nnf}. The relation between the finite presentability
of a coherent categorical structure and the finite presentability of
the monoid or group that arises from it remains unclear and much
remains to be done in this direction. 

\bibliographystyle{alpha}
\bibliography{papers}

\end{document}